\documentclass[final,review,authoryear,11pt]{elsarticle}
\usepackage{lipsum}
\makeatletter
\def\ps@pprintTitle{
 \let\@oddhead\@empty
 \let\@evenhead\@empty
 \def\@oddfoot{}%
 \let\@evenfoot\@oddfoot}
\makeatother

\usepackage{natbib}
\usepackage{amssymb}
\usepackage{multirow}
\usepackage{booktabs}
\usepackage{fullpage}
\usepackage{lscape}
\usepackage{tikz}
\usepackage{graphicx} 
\usepackage{amsmath}
\usepackage{amsthm}
\usepackage{float}
\usepackage{rotating}
\usepackage{todonotes}
\usepackage{algpseudocode}
\usepackage{algorithm2e}
\usepackage{graphicx}
\usepackage{caption}
\usepackage{subcaption}

\usepackage{xcolor,colortbl}
\usepackage{url}
\usepackage{enumitem}
\usepackage{xcolor}
\usepackage{tabu}

\begin{document}

\begin{frontmatter}



\title{{\LARGE Observational Data-Driven Modeling and Optimization\\ of Manufacturing Processes}}

\author {Najibesadat Sadati}
\ead{n.sadati@wayne.edu}
\author {Ratna Babu Chinnam$^{*}$}
\ead{Ratna.Chinnam@wayne.edu}
\author {Milad Zafar Nezhad.}
\ead{m.zafarnezhad@wayne.edu}
\address {Department of Industrial and Systems Engineering, Wayne State University, \\4815 Fourth Street, Detroit, MI 48202, USA}

\cortext[mycorrespondingauthor]{Corresponding author}

\begin{abstract}
The dramatic increase of observational data across industries provides unparalleled opportunities for data-driven decision making and management, including the manufacturing industry. In the context of production, data-driven approaches can exploit observational data to model, control and improve process performance. When supplied by observational data with adequate coverage to inform the true process performance dynamics, they can overcome the cost associated with intrusive controlled designed experiments and can be applied for both process monitoring and  improvement. We propose a novel integrated approach that uses observational data for identifying significant control variables while simultaneously facilitating process parameter design. We evaluate our method using data from synthetic experiments and also apply it to a real-world case setting from a tire manufacturing company.
\end{abstract}

\begin{keyword}
Parameter design \sep observational data \sep variable selection \sep data-driven modeling \sep response surface method \sep meta-heuristic optimization
\end{keyword}

\end{frontmatter}

\section{Introduction}
Parameter design is a methodology that seeks to improve the quality of products and processes. It categorizes variables affecting the desired performance characteristic into control variables and noise (uncontrollable) variables and sets the levels of the control variables so as to optimize the performance characteristic while minimizing the variance imposed by the noise variables \citep{robinson2004robust}. Among several methods developed in recent decades, response surface methodology (RSM) is seen to be an effective method for parameter design \citep{montgomery2008design}. In general, RSM methods rely on data from structured designed experiments for parameter design. These controlled experiments lead to greater confidence in the parameter design process. The drawback, however, is that these controlled experiments can be disruptive (e.g., to production during process parameter design) and cost prohibitive. 

On the contrary, most companies today already collect vast troves of process data but typically use them only for tracking purposes, not as a basis for improving operations \citep{Auschitzky2014}. For example, it is common in biopharmaceutical production flows to monitor more than 200 variables to ensure the purity of ingredients as well as the substances being made stay in compliance. One of the many factors that makes biopharmaceutical production and some other continuous processing industries so challenging is that yields can vary from 50\% to 100\% for no immediately discernible reason \citep{Columbus2014}. Using advanced analytics, a manufacturer was able to track the nine parameters that most explained yield variation and increase the vaccine’s yield by 50\%, worth between \$5M to \$10M in yearly savings for the single vaccine alone \citep{Auschitzky2014}. Initiatives such as Industrie 4.0, a German government initiative that promotes automation of the manufacturing industry with the goal of developing Smart Factories, and the rapid adoption of IoT (Internet of Things) devices in manufacturing can indeed drive innovation, competitiveness, and growth \citep{Löffler2013, Yosefi2014}. Although technology is revolutionizing the way in which data is generated and captured, there is need for analytics methods that can harness this data for enhanced monitoring, control, and management. We posit that ``observational'' data (i.e., data collected from routine production) often readily available from many modern production systems can facilitate adaptive process parameter design. This is particularly true if the process was historically operated under a variety of operating conditions due to existing control policies and/or manual interventions. The additional motivation for this research stems from a request by a leading global tire manufacturing company seeking assistance in improving the consistency of rubber compounds in terms of material viscosity.

Relying on historical observational data for parameter design brings about a number of challenges however it creates opportunity for reducing the cost of experimental design and eliminates process disruptions necessary for experiment execution to achieve a robust parameter set. This study seeks to develop an integrated approach for process parameter design based on observational data that jointly addresses the tasks of identifying critical controllable variables while simultaneously facilitating process parameter design. Our approach uses an iterative process that first determines the potential set of variables to control using a meta-heuristic optimization approach and then employs these candidate variables to build an appropriate RSM for the target variable of interest and finally optimizes the candidate control variable set for process performance robustness. The process is repeated till convergence. Our proposed method is a general approach and can be applied to different processes (systems) where the goal is to identify the most critical control variables and determine the appropriate settings for these control variables to optimize an output of interest (a quality characteristic) while minimizing the variation affected by uncontrollable noise variables. Our method is most applicable to manufacturing processes in different industries such as automotive, biotechnology, chemical and pharmaceutical industry. We demonstrate the performance of our method using data from synthetic experiments and discuss a real-world case study. 

The rest of this paper is organized as follows. Section II reviews the related literature regarding variable selection methods and robust parameter design. Section III explains our proposed approach for parameter design using observational data. Section IV describes validation results from synthetic experiments as well as a real-world case study. Finally, section V provides some closing remarks and identifies directions for future research.

\section{Literature Review}
Since this research seeks to develop an integrated approach for critical variable selection and robust parameter design using response surface methodology (RSM) applied to observational data, we first review some related literature. 

\subsection{Overview of Variable Selection Methods}
In the domains of statistics, data mining and machine learning, there is an extensive body of literature on the topic of variable or feature selection (in the rest of this section, we use the terms `feature' and `variable' interchangeably). These methods seek to achieve dimensionality reduction and improve model performance when the dataset carries several noisy and irrelevant variables. Depending on the nature of the modeling task, variable selection algorithms can be categorized as supervised \citep{weston2003use,song2007supervised}, unsupervised \citep{mitra2002unsupervised, dy2004feature}, and semi-supervised \citep{zhao2007semi, xu2010discriminative}. For a good review of these methods, see  \cite{jain198239,dash1997feature,liu2005toward,nakariyakul2009improvement}. Given our interest in parameter design, we limit our discussion to supervised variable selection methods. In a recent study, \cite{tang2014feature} reviewed many supervised variable selection algorithms/methods and classified them into: 1) Algorithms for `flat' features which assume that features are (somewhat) independent and most relevant for our study -- usually divided into three sub-groups: `filters', `embedded' methods, and `wrappers'; 2) Algorithms for `structured' features -- related to settings where features constitute/posses an intrinsic structure (e.g., group, tree, or graph structure); and 3) Algorithms for `streaming' features -- relevant to settings where the knowledge about the full feature space is unknown or dynamic (e.g., on-line applications). 

Feature selection methods can be compared based on several criteria. For example, filter methods rank features based on specific criteria such as Fisher score \citep{duda2012pattern} and mutual information (MI) \citep{kira1992practical, robnik2003theoretical}. Filters are computationally efficient and select features independently. The major disadvantage of filters is that they do not consider the impact of selected feature subset on the induction algorithm performance \citep{tang2014feature}. Wrapper models explicitly account for performance of the intended induction algorithm and search for the best variable subsets, yielding superior performance \citep{tang2014feature} but are iterative and tedious \citep{kohavi1997wrappers, inza2004filter}. Embedded methods accomplish variable selection as a part of the learning process itself \citep{guyon2003introduction}. While they combine the advantages of wrappers and filters, they are often limited to certain classes of models. Common ones are regularization regression models such as lasso, ridge and elastic-net regression, which jointly minimize fitting errors while also penalizing model complexity \citep{tang2014feature, nezhad2016safs}.

Many `hybrid' variable selection methods are also developed to seek synergies across filters, embedded methods, and wrappers. Most of these methods employ two phases. In the first phase they employ filter or embedded methods for ranking and reducing the number of variables and then, in the second phase, wrapper method is applied to select the desired number of variables among the reduced set of variables \citep{dash1997feature, uncu2007novel}. Here we present a brief review of hybrid methods due to their strong relevance to the proposed framework for parameter design. \cite{raymer2000dimensionality} proposed an integrated framework for feature selection, extraction, and classifier training using genetic algorithms (a ``derivative free" meta-heuristic optimization method). They employed a genetic algorithm to optimize `weights' assigned to features, which are used to rank the individual features. Das \cite{das2001filters} presented a hybrid algorithm that used boosted decision stumps as weak learners. Their method incorporates some of the benefits of wrappers, such as a natural stopping criterion, into a fast filter method. \cite{huang2007hybrid} developed a hybrid feature selection method using a genetic algorithm for finding a subset of variables which are most relevant to the classification task. Their method includes two optimization steps, the outer loop carries out `global' search to find best subset of features in a wrapper manner, and the inner loop performs a `local' search like a filter to improve the conditional MI, considered as an indicator for feature ranking. The results based on a real-world dataset demonstrate that this hybrid method is much more efficient than wrapper methods and outperforms filter methods in terms of accuracy. \cite{hsu2011hybrid} proposed a three-step framework that starts with preliminary screening and then continue with combination, and fine tuning. They applied their method to two bioinformatics problems and achieved better results in comparison to other methods. \cite{unler2011mr} presented a hybrid filter-wrapper variable selection method that relies on particle swarm optimization (PSO), another popular meta-heuristic optimization method, for support vector machine (SVM) classification. The filter method is based on MI and the wrapper method employs PSO for search. They applied the approach to some well-known benchmarking datasets with excellent accuracy and computational performance. 

In summary, hybrid approaches provide the best trade-off between computational burden and modeling accuracy in comparison with standard wrapper, filter and embedded methods. Meta-heuristic optimization methods (i.e., GA and PSO) seem to be quite effective for hybrid techniques. Given our goal to develop an integrated approach that not only selects significant variables for control but also executes parameter design, without loss of generality, we employ a hybrid approach for critical control variable selection. In the next section, we review some variable selection methods applied specifically for manufacturing process monitoring and parameter design. 

\subsection{Variable Selection for Manufacturing Process Monitoring and Parameter Design}
Variable selection in the context of manufacturing process monitoring has been studied by several researchers. Gauchi and Chagnon \cite{gauchi2001comparison} compared more than twenty variable selection methods in the context of partial least square (PLS) regression based on different real-world datasets related to chemical manufacturing processes. Among all methods, a stepwise variable selection based on maximum $Q_{cum}^{2}$ criterion outperformed the other methods. \textcolor{black}{In the other research, \cite{del2007system} present a simple, fast and network based approach for monitoring and parameter tuning of a high performance drilling process. They performed model validation using error-based performance indices and correlation analyses.} \cite{wang2009high} proposed a variable selection scheme based on a multivariate statistical process control (SPC) framework for fault diagnosis and process monitoring applications. The proposed framework selects the out of control variable first using a forward-selection algorithm and then monitors suspicious variables by setting up multi-variable charts. They demonstrated the effectiveness of their procedure based on a simulation study and a real-word experiment. In similar research, \cite{gonzalez2010variable} developed a two-state procedure for variable selection to monitor processes based on SPC. In the first stage, based on some criterion, all variables get sorted and in the second stage, variable selection is accomplished by calculating two measures (first measure is amount of residual information in unselected variables and the other measure is performance of control chart). Authors evaluated their method's effectiveness in a metal forming application and simulation study. \textcolor{black}{\cite{penedo2012hybrid} developed a new hybrid incremental model for manufacturing process monitoring. Their approach includes two iterative steps: a global model is first built using least squares regression and then employs a local model based on a fuzzy KNN smoothing algorithm. They implemented this hybrid model for tool wear monitoring in turning processes and the comparative results indicate computational efficiency and effectiveness.} More recently,  \cite{shao2013feature} presented a novel method for variable selection and parameter tuning in quality control of manufacturing processes. They used a cross-validation approach and consider false-positive and false-negative errors to identify the best subset of variables. They applied their method to data from an ultrasonic metal welding process of batteries achieving good monitoring performance. However, none of these methods directly investigate the relationship between control and noise variables in the context of process parameter design, the focus of this paper. 

\subsection{Parameter Design using Response Surface Methodology}
In the 1980s, \cite{taguchi1986introduction} proposed parameter design (a.k.a. robust design), which relies on statistical design of experiments for quality engineering. Taguchi recommended the reduction of data resulting from the experiments into signal-to-noise ratio measures for robust design. While this makes the approach relatively simple and practical, the statistical community embraced Response Surface Methodology (RSM) based methods and others for modeling the performance characteristic from the experimental data and then subjecting the resulting response function to more vigorous optimization for seeking robust designs and processes. \cite{montgomery2008design} describes RSM as “a collection of mathematical and statistical techniques useful for the modeling and analysis of problems in which a response of interest is influenced by several variables and the objective is to optimize this response”. These methods have been employed extensively in numerous industries. For a good review of RSM and its application, see \cite{myers2004response}. Two significant RSM approaches for robust design are: combined array designs and dual response surface approach. Both methods have their own advantages. For illustrative purposes, we employ the dual response surface approach, which is outlined in the next section as a part of the proposed approach.

\section{Proposed Integrated Variable Selection \& Parameter Design Approach}

\begin{figure}[H]
	\centering
	\includegraphics[trim={5cm 4.6cm 5cm 3.5cm},clip,scale= 0.7]{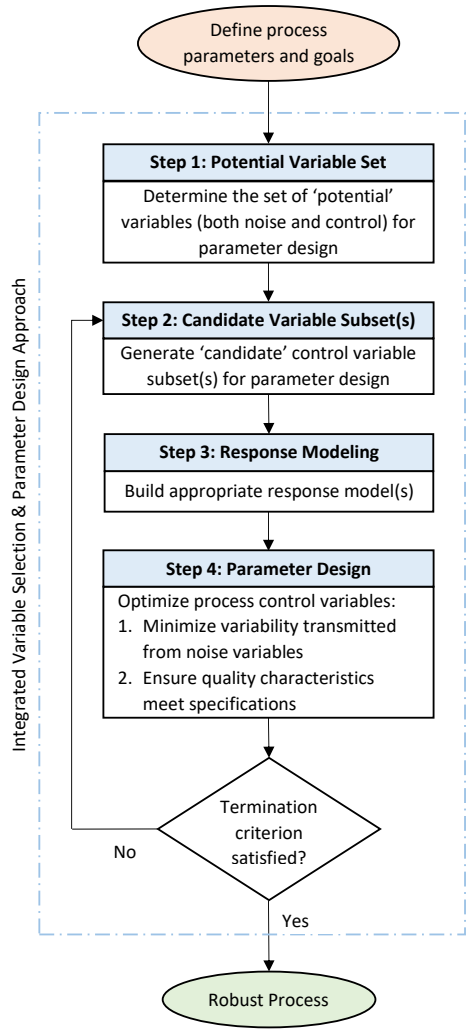} 
	\caption{Proposed integrated variable selection and parameter design approach}
	\label{figure1}
\end{figure}

This section outlines our proposed method for integrated critical variable selection and process parameter design using observational data. Figure \ref{figure1} depicts how these two tasks (i.e., variable selection and process parameter design) are integrated through a four-step recursive process: 
\begin{itemize}[noitemsep,topsep=0pt,parsep=0pt,partopsep=0pt]
	\item Step$-$1: Determine the set of `potential' variables for parameter design; 
	\item Step$-$2: Select or update `candidate' variable subset(s) for process control; 
	\item Step$-$3: Build appropriate response surface model(s) for the target performance characteristic of interest using the candidate variable subset(s); and 
	\item Step$-$4: Optimize the candidate process control variable subset(s) for robustness. 
Repeat steps 2$-$4 until the termination criterion is satisfied. Output the best parameter design settings.
\end{itemize} 

\subsection{Step$-$1: Determine set of `potential' control and noise variables of interest}
The proposed approach starts with identification of `potential' variables for consideration for robust/parameter design (both control and noise variables). This starting set is likely to contain some redundant variables that need elimination through the remaining steps of this approach. It is understood that historical data associated with these variables is readily available from the target process, collected at appropriate and meaningful intervals (e.g., for select production batches or sampled at reasonable intervals if the process is continuous) using standard manufacturing execution systems (MES). Variable selection methods such as filters can be used during this stage to roughly rank and select the list of variables affecting the target performance characteristic of interest (i.e., response variable). The variable set should be finalized using expert knowledge (e.g., using available process domain experts and production staff). The efficiency of this step greatly dictates the efficiency of the remaining steps of the overall approach. 

\subsection{Step$-$2: Determine `candidate' variable subset(s) for parameter design}
This step initiates the search for the subset of most effective control variables within the `potential' set from Step-1 for parameter design. Given the redundancy in the potential set, it constructs solution `candidate' subset(s) of controllable variables for parameter design. Both exact and meta-heuristic optimization methods can be employed for the search (i.e., construction and exploration of alternate variable subset(s)). \textcolor{black}{The proposed method is generally agnostic to the type of optimization technique employed. While exact methods can be computationally efficient, they might constrain the response surface models of Step-3 to linear or quadratic functions. Meta-heuristic optimization methods such as Genetic Algorithms (GA), Simulated Annealing, Ant Colony, Cross-Entropy Method and Swarm Optimization are ``derivative free" optimization methods and can support any arbitrary response surface modeling method for process optimization \citep{beruvides2016multi}}, albeit at the cost of increased computational burden. For example GA is a `parallel' search method and entertains multiple candidate variable subsets or solutions in any given iteration of Figure \ref{figure1}.

\subsection{Step$-$3: Build appropriate response surface model }
The proposed method is generic and can support any arbitrary method for modeling the process response mean and variance. Without loss of generality, we discuss here the dual response approach \citep{myers1973response} for modeling the process response. \textcolor{black}{The dual response surface methodology is a well-known approach in parameter design that considers two response surfaces, one for process mean and another for process variance, and then formulates them as primary and secondary response surfaces. The dual RSM needs less experiment runs and considers important control-by-control and control-by-noise interactions, which make it superior to classical Taguchi parameter design methods. The dual RSM formulates the regression model and constraint optimization to achieve robust parameter design \citep{vining1990combining,shaibu2009another}.} Without loss of generality, we illustrate the process using the following response model structure:
\begin{align}
&Y(\textbf{x},\textbf{z})= g(\textbf{x}) + f(\textbf{x},\textbf{z}) + \varepsilon.&
\end{align}

\noindent
where, $g(\textbf{x})$ is the function including just the control variables within the candidate set $\textbf{x}$ and $f(\textbf{x}, \textbf{z})$ is the function capturing the affects of noise variables $\textbf{z}$ as well as interactions between control and noise variables. \textcolor{black}{Assuming that the noise variables have a zero mean (can be standardized if necessary) the process mean response surface becomes:}
\begin{align}
&E_{\textbf{z}}\left[y(\textbf{x},\textbf{z})\right] = g(\textbf{x}).&
\end{align}

A method for obtaining variance response surface is to employ the transmission-of-error approach (so-called delta method) developed by \cite{myers2004response}. Per the delta method, one can derive the following variance response surface function: 
\begin{align}
&Var_{\textbf{z}}\left[y(\textbf{x},\textbf{z})\right]  =  \sigma_{\textbf{z}}^{2} \sum_{i=1}^{r} \Big[\dfrac{\sigma_{y(\textbf{x},\textbf{z})}}{\sigma_{z_{i}}}\Big]^{2} + \sigma^{2}&
\end{align}

\noindent
where, $\sigma_{\textbf{z}}^{2}$ is the variance of noise variables and $\sigma^{2}$ is the variance of residuals. 

For illustrative purposes, let us suppose that Eqn. (1) can be fit using the observational data through a quadratic regression model that consists of all first-order and second-order terms for the control variables, first-order terms for noise variables, and first-order interactions between noise and control variables. This regression model can be shown as follows: 
\begin{align}
&Y(\textbf{x},\textbf{z}) = \beta_{0} + \mathbf{x}^{'} \boldsymbol{\beta}_{1} + \textbf{x}^{'}\boldsymbol{\beta}_{2} \textbf{x}+\textbf{z}^{'}\boldsymbol{\beta}_{3}+ \textbf{x}^{'}\boldsymbol{\beta}_{4}\textbf{z}+\varepsilon.&
\end{align}      
\textcolor{black}{All parameters ($\beta_{i}$) in Eqn.(4) are obtained by response surface methodology.} Since for many applications the noise by noise interactions and second order terms of noise variables are not significant, we do not consider them in this illustrative response model.

\begin{table}[H]
	\scriptsize
	\centering
	\caption{Notation for process response regression model}
	\label{table1}
	
\begin{tabular}{ >{\small}c >{\small}l}
	\hline
	\textbf{Notation} &\textbf{Description} \\
	\hline
	$Y(\textbf{x},\textbf{z})$	& response variable\\
	$\textbf{x}$ & vector of control variables\\
	$\textbf{z}$ & vector of noise variables\\
	$\beta_{0}$ & model intercept \\
	$\boldsymbol{\beta}_{1}$ & coefficient vector for $1^{st}$ order control variable terms\\
	$\boldsymbol{\beta}_{2}$ & coefficient matrix for $2^{nd}$ order control variables terms and control by control interactions\\
	$\boldsymbol{\beta}_{3}$ & coefficient vector for $1^{st}$ order noise variable terms\\
	$\boldsymbol{\beta}_{4}$ & coefficient matrix for control by noise interactions\\
	$\varepsilon$ &error assumed to be an $i.i.d.$ random variable with $N (0, \sigma^{2})$ distribution\\
	\hline	
\end{tabular}

\end{table}

\textcolor{black}{Given the quadratic regression model in (4), the resulting response surface models for the process mean and variance are as follows:}
\begin{align}
&E_{\textbf{z}}\left[y(\textbf{x},\textbf{z})\right]  = \beta_{0} +\textbf{x}^{'}\boldsymbol{\beta}_{1} + x^{'}\boldsymbol{\beta}_{2} x.&
\end{align}
\textcolor{black}{
\vspace{-10mm}
\begin{align}
&Var_{\textbf{z}}\left[y(\textbf{x},\textbf{z})\right]= (\boldsymbol{\beta}_{3}^{'}+\textbf{x}^{'}\boldsymbol{\beta}_{4})Var_{\textbf{z}}[z] (\boldsymbol{\beta}_{3}+\boldsymbol{\beta}_{4}^{'} \textbf{x})+ \sigma^{2}&
\end{align}    
\noindent
where 〖$Var_{\textbf{z}}[z]$ denotes the variance-covariance matrix of $z$. If we suppose that $Var_{\textbf{z}}[z]= \sigma_{\textbf{z}} ^{2}I$, then Eqn. (6) can be simplified as:
\begin{align}
&Var_{\textbf{z}}\left[y(\textbf{x},\textbf{z})\right]= \sigma_{\textbf{z}}^{2} (\boldsymbol{\beta}_{3}^{'}+\textbf{x}^{'}\boldsymbol{\beta}_{4}) (\boldsymbol{\beta}_{3}+\boldsymbol{\beta}_{4}^{'} \textbf{x})+ \sigma^{2}.&
\end{align}}\label{eq6}

\vspace{-10mm}
\subsection{Step$-$4: Parameter design}
Our task here is to optimize the control variables using the response surface models from Step$-$3 to achieve robustness in process performance. In particular, we should apply constrained optimization to find optimal process parameters that ensure that the process response mean is as close to the desired target ($t$) as possible (forms a `constraint') while minimizing the process variance (forms the `objective' function). One way to setup the formulation for this optimization problem is as follows:

\textit{Parameter Design Formulation \#1}

$Min$ \quad $Var_{\textbf{z}}\left[y(\textbf{x},\textbf{z})\right]$ 

$Subject\:to$ \quad $E_{\textbf{z}}\left[y(\textbf{x},\textbf{z})\right]  = t$

\qquad \qquad \qquad  $\textbf{x} \in R^{k}$ 

An alternative formulation would seek to balance the deviation of the process mean from the target against process response variance. This is also necessary if there exist no control parameter settings that will yield the target response. 

\textit{Parameter Design Formulation \#2}

$Min$ \quad $\alpha \cdot Var_{\textbf{z}}\left[y(\textbf{x},\textbf{z})\right] + (1-\alpha) \cdot (E_{\textbf{z}}\left[y(\textbf{x},\textbf{z})\right]-t)^{2} $\

$Subject\:to \quad \textbf{x} \in R^{k}$ and $\alpha \in (0,1) $

\noindent
where $\alpha$ controls the importance assigned to variance minimization vs. deviation of expected process response from the target. \textcolor{black}{The main drawback of this method is the difficulty in selecting an appropriate value for $\alpha$. There exist several approaches in the literature under the domain of multi-objective optimization to identify the appropriate weights ($\alpha$) for the objective, such as rating methods, ranking methods, categorization methods, and ratio questioning \citep{marler2010weighted}.}

Another alternative would be to employ the signal-to-noise (SNR) ratios (see \cite{taguchi1986introduction} and \cite{Murphy2008} for more details) as objectives for optimization. Depending on the target response characteristic, they can take the following forms \citep{box1988signal}:

\begin{description}[noitemsep,topsep=0pt,parsep=0pt,partopsep=0pt]
\item\small $\bullet$ SNR$_{T}$: ``Nominal is best" setting -- goal is to achieve a specific target value for the response variable while minimizing the variability.
\begin{align}
	&SNR_{T}=10 \log_{10} \Big(\frac{\bar{y}^2}{s^2}\Big)&
 \end{align} 
where $\bar{y}$ = $\frac{1}{n}$ $\sum_{i=1}^{n} {y_i}$, $s = (\frac{1}{n-1} \sum_{i=1}^{n} {(y_i-t})^2)^{1/2}$, ${y_i}$ denotes the $i^{th}$ response/observation, and $n$ denotes number of observations.

\item\small $\bullet$ SNR$_{L}$: ``Larger is better" setting -- goal is to maximize the response variable (e.g., durability) while minimizing the variability.
\begin{align}
	&SNR_{L}=-10 \log_{10}  \Big[\frac{1}{n} \sum_{i=1}^{n} {\Big(\frac{1}{y_i}}\Big)^2 \Big].&
\end{align} 

\item\small $\bullet$ SNR$_{S}$: ``Smaller is better" setting -- goal is to minimize the response variable (e.g., impurity) while minimizing the variability.
\begin{align}
	&SNR_{S}=-10 \log_{10}  \Big[\frac{1}{n} \sum_{i=1}^{n} {y_i^2}\Big].&
\end{align} 
\end{description}

\noindent See Nair \cite{nair1992taguchi} for detailed discussion on the pros and cons of employing different approaches (SNRs vs RSM methods and more) for parameter design.

Without loss of generality, to simplify the illustration of the framework, we consider the setting with a single process response variable of interest. The proposed method can be generalized for the case of multiple responses. One approach would be to formulate the multi-response problem as a constrained optimization problem, where one response forms the objective function and the others are handled through constraints. See \cite{antony2001simultaneous} for an example application. Goal programming and other general purpose multi-objective optimization methods can also be employed. See \cite{deb2016multi} for a good overview of these methods.

\subsection{Termination}

If the parameter settings resulting from Step$-$4 for the (best) candidate subset are not satisfactory, we store these results (candidate variable subset(s) and resulting performance(s)) for future recall and go back to Step$-$2 for generating alternate candidate variable subset(s) from the potential variable set pool of Step$-$1 and repeat Steps 3 and 4. The process is repeated until the termination criterion is satisfied (can be a threshold for acceptable process variance or SNR or predetermined number of iterations) or if there is no improvement in the quality of the parameter design setting over multiple iterations. 

The pseudo code for our approach is provided below:

\begin{algorithm}[H]
	\caption{Pseudo code for integrated variable selection and parameter design approach}
	\KwData{Supply observational dataset, Sets of candidate noise variables $(\mathbf{Z})$ and control variables $(\textbf{X})$, Goals for parameter design}
	\KwResult{Robust process parameter design setting $\hat{\textbf{x}}$}
	\While{Termination criterion is not met}{
		Select subset(s) of control variables for parameter design $(\underline{\textbf{X}})$ using a meta-heuristic optimization method\;
		Build appropriate response surface model $\hat{y}(\underline{\textbf{X}},\textbf{Z})$\;               
		Conduct parameter design to identify best control variable settings $\hat{\textbf{x}}$\;
		Evaluate fitness of $\hat{\textbf{x}}$ for satisfactory performance\;
	}
\end{algorithm}

\section{Model Implementation \& Validation}
To validate the effectiveness of the proposed approach, we first test the performance of the proposed approach on datasets derived from synthetic experiments and compare its performance with some traditional alternative approaches. Later, we also report results from a real-world case study that stems from a tire manufacturing company.

\subsection{Synthetic Experiments}

In this section, we test the performance of the proposed method using data from synthetic experiments. We consider different number of control, noise, and `dummy' variables and we present numerical comparison of the results from our approach to two popular baseline methods: 1) filter based variable selection with mutual information (MI) criterion; and 2) random forest method. In particular, we evaluate 18 different experimental instances of growing complexity, with each test setting being replicated twenty times for evaluating robustness of the approach.

\clearpage
\begin{landscape}
\centering 
\begin {table}[H]
\caption {Experimental settings for synthetic validation study}\label{table2} 
\begin{tabular}{ c c c c c c c c c c c c c c c c c | c c }
	\hline
	\multicolumn{17}{c|}{Experimental Setting}  &\multicolumn{2}{c}{GA Setting}\\	
	$\# Exp.$ & $\# C$ & $\# N$ & $\# D$ & $\# O$ & $\sigma_{\varepsilon}$ & $\sigma_{\textbf{X}}$ & $E(\textbf{X})$ & $\sigma_{\textbf{Z}}$ & $E(\textbf{Z})$ & $\sigma_{\textbf{D}}$  & $E(\textbf{D})$ & $\sigma_{\boldsymbol{\beta}_1}$ & $\sigma_{\boldsymbol{\beta}_2}$ & $\sigma_{\boldsymbol{\beta}_3}$& $\sigma_{\boldsymbol{\beta}_4}$& $\# Runs$ & Pop. & Gens. \\
	\hline
	1 & 2 & 2 & 2 & 100 & 1 & 2 & 0 & 1 & 0 & 2 & 0 & 0.1 & 0.5 & 2 & 0.5 & 20 & 4 & 10\\ 
	2 & 2 & 2 & 2 & 100 & 1 & 2 & 1 & 1 & 0 & 2 & 1 & 0.1 & 0.1 & 2 & 0.1 & 20 & 4 & 10\\ 
	3 & 2 & 2 & 2 & 100 & 1 & 2 & 10 & 1 & 0 & 2 & 10 & 1 & 0.01 & 2 & 0.01 & 20 & 4 & 10\\ 
	\hline
	4 & 2 & 2 & 2 & 100 & 1 & 8 & 0 & 1 & 0 & 8 & 0 & 0.1 & 0.5 & 8 & 0.5 & 20 & 4 & 10\\ 
	5 & 2 & 2 & 2 & 100 & 1 & 8 & 1 & 1 & 0 & 8 & 1 & 0.1 & 0.1 & 8 & 0.1 & 20 & 4 & 10\\ 
	6 & 2 & 2 & 2 & 100 & 1 & 8 & 10 & 1 & 0 & 8 & 10 & 1 & 0.01 & 8 & 0.01 & 20 & 4 & 10\\ 
	\hline
	7 & 4 & 4 & 4 & 1000 & 1 & 2 & 0 & 1 & 0 & 2 & 0 & 0.1 & 0.5 & 2 & 0.5 & 20 & 10 & 25\\ 
	8 & 4 & 4 & 4 & 1000 & 1 & 2 & 1 & 1 & 0 & 2 & 1 & 0.1 & 0.1 & 2 & 0.1 & 20 & 10 & 25\\ 
	9 & 4 & 4 & 4 & 1000 & 1 & 2 & 10 & 1 & 0 & 2 & 10 & 1 & 0.01 & 2 & 0.01 & 20 & 10 & 25\\
	\hline
	10 & 4 & 4 & 4 & 1000 & 1 & 8 & 0 & 1 & 0 & 8 & 0 & 0.1 & 0.5 & 8 & 0.5 & 20 & 10 & 25\\ 
	11 & 4 & 4 & 4 & 1000 & 1 & 8 & 1 & 1 & 0 & 8 & 1 & 0.1 & 0.1 & 8 & 0.1 & 20 & 10 & 25\\ 
	12 & 4 & 4 & 4 & 1000 & 1 & 8 & 10 & 1 & 0 & 8 & 10 & 1 & 0.01 & 8 & 0.01 & 20 & 10 & 25\\
	\hline
	13 & 8 & 8 & 8 & 1000 & 1 & 2 & 0 & 1 & 0 & 2 & 0 & 0.1 & 0.5 & 2 & 0.5 & 20 & 20 & 60\\ 
	14 & 8 & 8 & 8 & 1000 & 1 & 2 & 1 & 1 & 0 & 2 & 1 & 0.1 & 0.1 & 2 & 0.1 & 20 & 20 & 60\\ 
	15 & 8 & 8 & 8 & 1000 & 1 & 2 & 10 & 1 & 0 &2 & 10 & 1 & 0.01 & 2 & 0.01 & 20 & 20 & 60\\
	\hline
	16 & 8 & 8 & 8 & 1000 & 1 & 8 & 0 & 1 & 0 & 8 & 0 & 0.1 & 0.5 & 8 & 0.5 & 20 & 20 & 60\\ 
	17 & 8 & 8 & 8 & 1000 & 1 & 8 & 1 & 1 & 0 & 8 & 1 & 0.1 & 0.1 & 8 & 0.1 & 20 & 20 & 60\\ 
	18 & 8 & 8 & 8 & 1000 & 1 & 8 & 10 & 1 & 0 &8 & 10 & 1 & 0.01 & 8 & 0.01 & 20 & 20 & 60\\
	\hline
\end{tabular}
\end {table}
\end{landscape}
\clearpage

Based on the number of variables, we have three groups of experiments: 1) two control and two noise variables, 2) four control and four noise variables, and 3) eight control and eight noise variables. The settings for the different instances are given in Table (2) in which $C$, $N$, and $D$ denote the number of control variables, noise variables, and pure dummy variables, respectively. Dummy variables do not have any impact on the process response and are pure white noise variables introduced to test the ability of the proposed approach in filtering out irrelevant variables. $O$ denotes the number of observations within the observational dataset. 

For these experiments, during each replication, we generate the model coefficients of variables and interactions for the `true' model (i.e., Eqn. (4)) as Gaussian distributed random values with zero mean. As reported in Table 2, in seeking  meaningful and revealing experiments, the scales for the variances of model coefficients are determined based on the number of control and noise variables for each experiment. For instance, the expected value for each of the control and dummy variables in experiment \#3 is 10.0, therefore, we considered a relatively low value of a variance of 0.01 for each of the quadratic coefficients ($\boldsymbol{\beta}_{2}$ and $\boldsymbol{\beta}_{4}$) to balance their effect on the response variable ($Y$). 

The experiments were conducted with the ground truth process response being driven by a quadratic response model as in Eqn. (4), with the goal of minimizing the process response variance while maintaining the expected response at the desired target (i.e., we employed Parameter Design Formulation \#1 for Step$-$4). There are multiple goals for these experiments. The first objective is to evaluate the effectiveness of the proposed method in its ability for control variable selection, i.e., the ability to differentiate the true control/noise variables from the pure dummy variables. The second goal is to evaluate the quality of the recommended parameter design in relation to the truly optimal setting derived from the ground truth model. The target response $t$ is identified for each instance as the expected response at the robust control parameter setting that leads to the least process variance using the ground truth model. 

For Step$-$2, we employed a standard Genetic Algorithm (GA), which is responsible for control variable selection. For more details regarding genetic algorithms, readers can refer to \cite{whitley1994genetic}. The parameters for the GA algorithm are also listed in Table 2 to the far right (Pop. denotes size of the candidate subset population and Gens. denotes the number of generations employed for the search). For Step$-$3, we employed the dual RSM approach. Finally, we have applied three criteria to investigate the performance of the proposed approach: 1) Recognition ratio for true control variables, 2) False-recognition ratio of non-relevant dummy variables as relevant variables, and 3) Accuracy of the RSM model coefficients in relation to the coefficients in the ground true model.

As for a baseline (`traditional') approach, we employ a typical sequential process where the set of `potential' variables are first reduced to a smaller subset of more relevant variables that will then be employed by the RSM and parameter design stages. We have employed two popular techniques for variable selection (i.e., reducing the potential set): 1) Filter method based on mutual information is first used to screen important control variables, and 2) Random forest is used for ranking the variable importance and then select a subset of high-quality control variables. We used the number of real control variables ($C$) as the threshold for variable selection among all ranked variables to minimize redundancy and to maximize relevant features and for giving these alternate methods full advantage (proposed method is not informed that there are $C$ real control variables).

\begin{figure}[H]
	\centering
	\begin{subfigure}{.5\textwidth}
		\centering
		\includegraphics[width=0.98\linewidth]{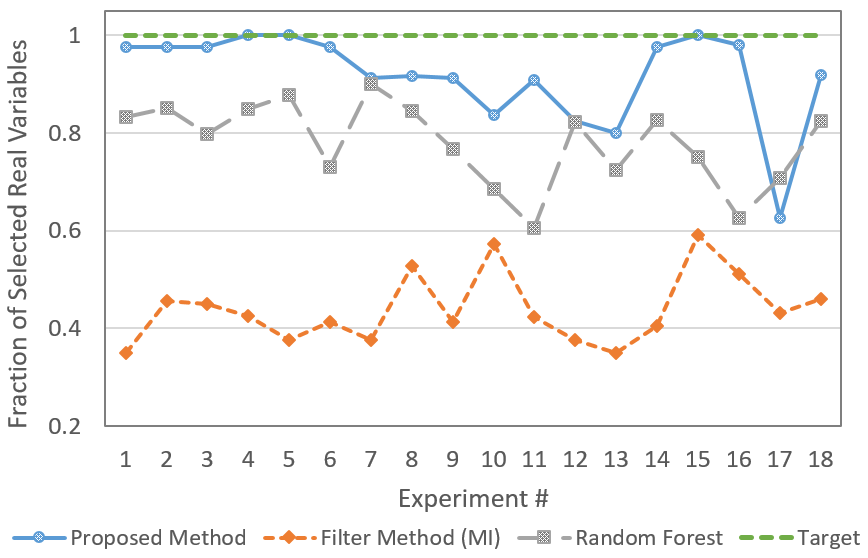}
		\caption{Fraction of Selected Real Variables (Target=1)}
		\label{fig:sub1}
	\end{subfigure}%
	\begin{subfigure}{.5\textwidth}
		\centering
		\includegraphics[width=0.98\linewidth]{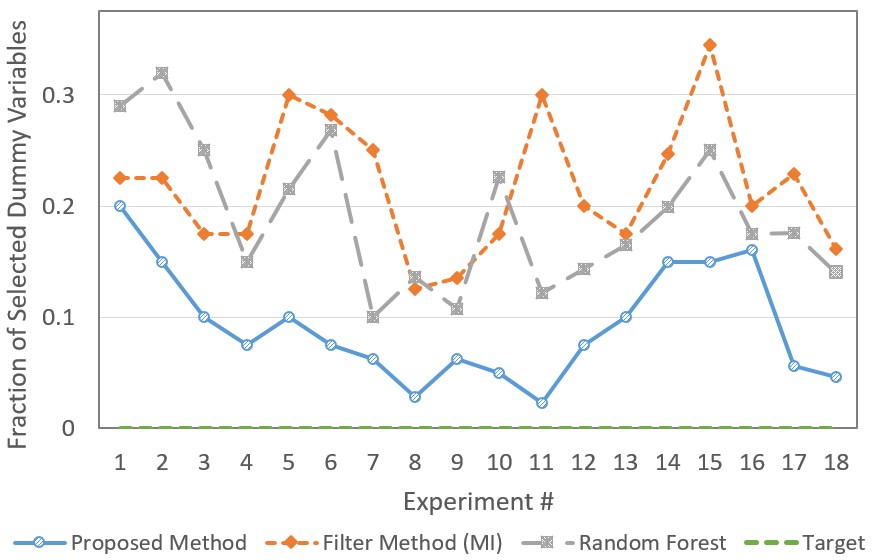}
		\caption{Fraction of Selected Dummy Variables (Target=0)}
		\label{fig:sub2}
	\end{subfigure}
	\caption{Comparison of proposed method with traditional sequential approaches}
	\label{fig:test}
\end{figure}

Figure \ref{fig:test} reports the performance results for both the proposed approach as well as the two baseline approaches. All results are based on 20 replications for each instance. Figure \ref{fig:test}(a) reports the recognition ratio for true control variables where as \ref{fig:test}(b) reports the false-recognition ratio for dummy variables. It is clear that the proposed approach significantly outperforms the baseline methods in almost all instances. 

\begin{figure}[H]
	\centering
	\includegraphics[height=3.6in,trim={1.3cm 0.8cm 1.3cm 0.9cm},clip]{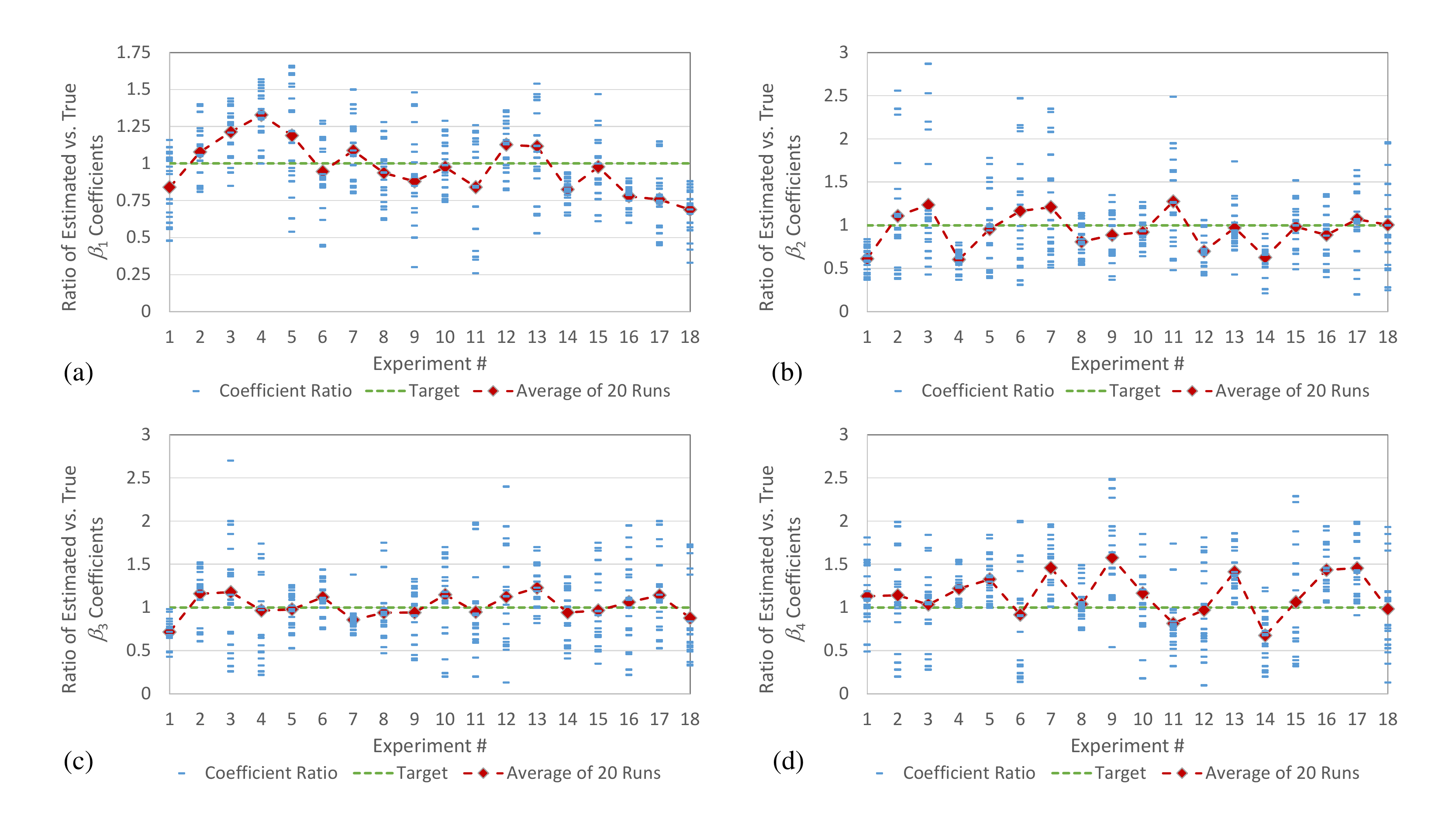}
	\caption{Mean Ratio of Estimated vs. True RSM Coefficients (Target = 1)}
	\label{figure3}
\end{figure}

To investigate the accuracy of the resulting RSM model, we compared estimated model parameters/coefficients with associated true parameters of the ground truth model. For brevity, we evaluate and report the proportional average of each generated parameter to corresponding parameter in the estimated model for each of the $\beta$ groups (i.e., we calculate $\beta_{i}/\hat{\beta_{i}}$ for each parameter within the coefficient group $i$ and average these ratios for each group). Figure \ref{figure3} reports modeling coefficient accuracy results for all 18 experimental instances. The average ratios across the 20 replication runs are quite close to the desired target ratio of 1.0 for all the 18 experiments with reasonable fluctuation in coefficient ratios. Overall, we can declare that the proposed integrated approach is rather effective in learning/identifying the real effects within the observational dataset to facilitate parameter design.

\subsection{Case Study}
We report here results from implementing our proposed approach on a real-world case study using observational data collected by a \textcolor{black}{pneumatic} tire manufacturing company. 

\textcolor{black}{Tires are manufactured according to relatively standardized processes and machinery by all the manufacturers. Tire is an assembly of numerous components that are typically built-up on a drum and then cured in a press under heat and pressure. Heat facilitates a polymerization reaction that crosslinks rubber monomers to create long elastic chains. Bulk raw materials such as synthetic and natural rubber, carbon black (gives reinforcement and abrasion resistance), and a variety of chemicals (e.g., silica, sulfur, vulcanizing accelerators, activators, and antioxidants) are employed to produce numerous specialized components (e.g., tread, side wall, inner liner) that are assembled and cured into a tire. The production process typically involves five stages: 1) Compounding/mixing operation brings together all the ingredients required to mix a batch of rubber compound, 2) Component preparation typically involves extrusion and bead building, 3) Tire building assembles the components onto a tire building drum, 4) Curing is the process of applying pressure and heat to stimulate the chemical reaction between the rubber and other materials, and 5) Final finish and testing.} 

\textcolor{black}{The case study is focused on a particular mixing line that produces a certain rubber compound for a specific tire component. As noted earlier, rubber compounds are produced by mixing a variety of chemicals in particular proportions along with rubber and fillers. The mixing process is typically a batch operation and involves several sub-stages. First stage involves mixing of carbon black along with rubber and other chemicals. Raw material properties, weights, material supplier location, and mixer equipment parameters, all influence the batch quality. The second stage involves addition of silica to the compound. The compound undergoes few additional stages  of processing before being ready for consumption to produce tire components. In particular, this case study focuses on the ``carbon mixing stage". The quality of the batch depends on the input material weights, properties, suppliers locations, and process parameters (e.g., chamber temperature, ram pressure, mix power, rotor speed, batch time and so on).}
	
\textcolor{black}{The most significant quality characteristic of interest for the carbon mixing stage is the rubber compound viscosity (a measure of its resistance to gradual deformation by shear or tensile stress, in other words, resistance to flow), which is necessary for effective control of subsequent processes as well as the performance of the tire. Of particular interest is the ``minimum viscosity" measurement from the viscosity test (which also yields other measures). Since uncontrollable factors such as ambient room temperature/moisture and raw material properties affect the compound's minimum viscosity, it is essential to minimize the impact of these factors on the batch's minimum viscosity (our target) by identifying the most important controllable variables and also determining the optimal settings for these variables to produce consistent batches with the desired level of minimum viscosity.} 

The dataset carries information regarding a variety of raw materials and process parameters (including input material weights and properties, raw material sources, mixing conditions and other parameters) and the goal is to control the minimum viscosity of the compound around a target value of 65 units. The dataset consisted of 214 regularly sampled observations collected over a span of several months. For confidentiality reasons, we are unable to disclose full details. 

\textcolor{black}{In order to apply the proposed approach to this production line, after data preprocessing, we employed a quadratic RSM that also accounts for the covariance terms (based on Eqn.(6)) for modeling the minimum viscosity. We also applied a standard genetic algorithm for key control variable selection and parameter design.} In particular, it involved executing the steps outlined in Algorithm 2.

\begin{algorithm}[t]
	\caption{Pseudo code of implementing proposed approach for tire plant}
	\KwData{Supply observational dataset, Sets of candidate noise variables $(\textbf{Z})$ and control variables $(\bar{\textbf{X}})$, Goals for parameter design}
	\KwResult{Robust process parameter design settings $\hat{x}$}
	\textbf{Initialize:} Consultation with process engineers revealed 6 noise variables and 10 `potential' control variables that could affect minimum viscosity of rubber compound.\\
	\While{Termination criterion is not met}{
		Choose random initial population or generate next generation of control variable subsets (i.e., `individuals' for the generation) using GA algorithm\;
		Build appropriate RSM model(s) using noise variables and chosen control variables for each individual of the generation\;               
		Evaluate fitness of each individual within the generation based on the parameter design optimization model in order to minimize variability transmission from noise variables while simultaneously ensuring that mean of minimum viscosity meets specification\;
	}
\end{algorithm}


The proposed method identified six of the ten potential control variables to be critical. For consistency, the two baseline methods (i.e., the Random Forest and the MI based Filter) were also asked to identify the six best control variables. The three methods yielded different sets of critical control variables. The Random Forest and Filter methods identified three and two variables that are in common with the proposed method, respectively. The results of the case study can be analyzed from two aspects. First, we pay attention to the quality characteristic of response value. As reported in Table \ref{table 3}, the proposed approach did improve the robustness of the process in maintaining the mean compound viscosity around the target value (65 units). Second, the variance of the response is much improved (\textcolor{black}{18.7}) in relation to the two baseline methods. Overall, the proposed integrated approach for control variable selection and parameter design seems to significantly outperform the traditional methods in terms of process quality and consistency for this real-life case study.

\begin{table}[H]
	\small
	\centering
	\taburulecolor{black}
	\caption{Process Result Comparison for Different Parameter Design Methods}
	\label{table 3}
	{\color{black}
	\begin{tabular}{|c | c |c| }
		\hline
		\textbf{Method} & \textbf{Mean Response} & \textbf{Variance of Response} \\ 
		\hline
		\multicolumn{1}{|c|}{\textit{Ideal Target}}& 65.0 & 0.0  \\
		\hline
		\multicolumn{1}{|c|}{Proposed Method}& 64.4 & 18.7  \\
		\hline
		\multicolumn{1}{|c|}{Filter Method with MI}  & 75.0 & 33.7   \\
		\hline
		\multicolumn{1}{|c|}{Random Forest}& 62.7 & 27.6\\ 
		\hline
	\end{tabular}}
\end{table} 

\section{Conclusion}
Manufacturing companies today are collecting vast troves of process data but typically use them only for monitoring purposes and not as a basis for improving operations. This data collection trend will accelerate as companies further embrace initiatives such as Industry 4.0 and adoption of IoT devices for innovation, competitiveness, and growth. \textcolor{black}{Extant literature has not addressed this gap to develop a framework that can use routinely collected observational data to optimize the process parameter design and improve the process performance. The main contribution of this research is the development of a novel integrated approach that: 1) employs readily available observational data from routine production to achieve robust process parameter design in manufacturing processes and 2) simultaneously discovers the more important control variables.} Traditional controlled designed experiments can be challenging in real-world production environments and this paves for an effective alternative approach to attain robust process parameter conditions. The proposed framework relies on an integrated control variable selection, response surface modeling, and optimization methodology. We report promising results from a synthetic experimental study as well as illustrative results from a tire compound production process case study.

There are several limitations to the proposed method. The quality of the results are limited by the span/coverage of the observational data with regard to the critical variables that truly impact the process performance. Future research can consider extending the proposed approach to include limited controlled experiments, in the lines of a hybrid evolutionary operation technique. Another stream of research can explore the path of robust optimization or stochastic programming rather than relying on meta-heuristic optimization methods for parameter design and variable selection.

\section*{References}
\small
\bibliographystyle{model5-names}\biboptions{authoryear, round}
\bibliography{References}{}

\end{document}